\documentclass{amsart}

\def\eps{\varepsilon }

\def\gam{\gamma }
\def\RR{\mathbb R}

\def\ZZ{\mathbb Z}

\newcommand{\set}[1]{\{#1\}}%set
\providecommand{\abs}[1]{\lvert#1\rvert}
\providecommand{\norm}[1]{\lVert#1\rVert}
\newcommand{\remove}[1]{ }

\newtheorem{theorem}{Theorem}%[section]

\newtheorem{lemma}[theorem]{Lemma}
\newtheorem{corollary}[theorem]{Corollary}

\theoremstyle{remark}
\newtheorem*{remark}{Remark}
\newtheorem*{remarks}{Remarks}

\numberwithin{equation}{section}

\begin{document}
\title[Discrete Ingham type inequalities]{Discrete Ingham type inequalities with a weakened gap condition}
\author{Vilmos Komornik}
\address{D\'epartement de math\'ematique\\
          Universit\'e Louis Pasteur\\
          7, rue Re\-n\'e Descartes\\
          67084 Strasbourg Cedex, France}
\email{komornik@math.u-strasbg.fr}
\author{Paola Loreti}
\address{Dipartimento di Metodi e Modelli\\
            Matematici per le Scienze Applicate\\
            Universit\`a degli Studi di Roma ``La Sapienza''\\
            Via A. Scarpa, 16\\
            00161 Roma, Italy}
\email{loreti@dmmm.uniroma1.it}
\keywords{Observability, Fourier series, vibrating strings}
\date{Version of February 2nd, 2007}

\begin{abstract}
We establish discrete Ingham type and Haraux type inequalities for exponential sums satisfying a weakened gap condition. 
They enable us to obtain discrete simultaneous observability theorems for systems of vibrating strings or beams.
\end{abstract}

\maketitle

\section{Introduction}\label{s1}

Harmonic analysis is an efficient tool in control theory. For example, various generalizations of a classical theorem of Ingham proved to be very helpful in establishing crucial observability theorems in many recent studies; see, e.g., \cite{KomLor}. The purpose of this work is to prove a new generalization, by weakening a usual gap condition and by replacing the observed integrals by Riemann sums, more realistic from a practical point of view. The theorem is then applied for the solution of two observability problems concerning systems of strings or beams.

Given a strictly increasing sequence $(\omega_k)_{-\infty}^{\infty}$ of real numbers, we consider functions of the form
\begin{equation}\label{11}
x(t)=\sum_{k=-\infty}^{\infty} x_ke^{i\omega_kt}
\end{equation}
with complex coefficients $x_k$. The following generalization of Parseval's equality, which improved an earlier result of Jaffard, Tucsnak and Zuazua \cite{JafTucZua1997}, was established in 
\cite{BaiKomLor103}:

\begin{theorem}\label{t1}
Assume that there exists a positive number $\gam$ satisfying
\begin{equation}\label{12}
\omega_{k+2}-\omega_k\ge2\gam\quad\text{for all}\quad k.
\end{equation}
Fix $0<\gam_0\le \gam$ arbitrarily and and set
\begin{align*}
&A_1:=\set{k\in\ZZ\ :\ \omega_{k}-\omega_{k-1}\ge\gam_0\quad\text{and}\quad \omega_{k+1}-\omega_k\ge\gam_0};\\
&A_2:=\set{k\in\ZZ\ :\ \omega_{k}-\omega_{k-1}\ge\gam_0\quad\text{and}\quad \omega_{k+1}-\omega_k<\gam_0}.
\end{align*}
Then for every bounded interval $I$ of length $|I|>2\pi/\gam $ there exist two positive constants $c_1$ and $c_2$ such 
that.
\begin{equation*}
c_1Q(x)\le\int_I \abs{x(t)}^2\ dt \le c_2Q(x)
\end{equation*}
for all sums of the form \eqref{11} with square summable coefficients,
where we use the notation
\begin{equation*}
Q(x):=\sum_{k\in A_1}\abs{x_k}^2
+\sum_{k\in A_2}\abs{x_k+x_{k+1}}^2+(\omega_{k+1}-\omega_k)^2(\abs{x_k}^2+\abs{x_{k+1}}^2).
\end{equation*}
\end{theorem}

\begin{remark}
Under the stronger assumption 
\begin{equation}\label{13}
\omega_{k+1}-\omega_k\ge\gam\quad\text{for all}\quad k
\end{equation}
our result reduces to a classical theorem of Ingham \cite{Ing}:
\begin{equation*}
c_1\sum_{k=-\infty}^{\infty}\abs{x_k}^2 \le \int_I \abs{x(t)}^2\ dt \le c_2\sum_{k=-\infty}^{\infty}\abs{x_k}^2.
\end{equation*}
\end{remark}

We shall establish the following discrete version of Theorem \ref{t1}:

\begin{theorem}\label{t2}
Assume that there exists a positive number $\gam$ satisfying \eqref{12}, and introduce the sets $A_1$, $A_2$ as above.
Given $0<\delta\le \pi/\gam$ arbitrarily, fix an integer
$J$ such that  $J\delta>\pi/\gam$. Then there exist two positive constants $c_1$ and $c_2$,  depending only on $\gam$ and $J\delta$, such that
\begin{equation}\label{14}
c_1Q(x)\le
\delta\sum_{j=-J}^J\abs{x(j\delta)}^2
 \le c_2Q(x)
\end{equation}
for all functions \eqref{11} whose coefficients satisfy the condition 
\begin{equation}\label{15}
x_k=0\quad\text{whenever}\quad \abs{\omega_k}>\frac{\pi}{\delta}-\frac{\gam}{2}.
\end{equation}
\end{theorem}

\begin{remarks}\mbox{}

\begin{itemize}
\item Under the stronger gap condition \eqref{13}, Theorem \ref{t2} reduces to an earlier result proved in 
\cite{KomLor132}:
\begin{equation*}
c_1\sum_{k=-\infty}^{\infty}\abs{x_k}^2 \le 
\delta\sum_{j=-J}^J\abs{x(j\delta)}^2
\le c_2\sum_{k=-\infty}^{\infty}\abs{x_k}^2
\end{equation*}
instead of \eqref{14}. 

\item In view of Remarks 2.6 and 2.7 in \cite{BaiKomLor103} \eqref{14} implies that, more generally,
\begin{equation*}
c_1Q(x)\le\delta\sum_{j=-J}^J\abs{x(t'+j\delta)}^2\le c_2Q(x)
\end{equation*}
for every $t'\in\RR$, for all functions \eqref{11} whose coefficients satisfy the condition \eqref{15}.
The constants $c_1$, $c_2$ depend only on $t'$, $\gam$ and $J\delta$.

\item Theorem \ref{t1} follows from Theorem \ref{t2}. Indeed, fix a bounded interval  $I=[t'-R,t'+R]$ with $R>\pi/\gam$, choose $\delta=R/J$ for every sufficiently large positive integer $J$, and let $J\to\infty$ in the resulting estimates.
\end{itemize}
\end{remarks}

In  the sequel we often write $A\asymp B$ instead of double inequalities of the form $c_1A\le B\le c_2A$ for brevity.

The preceding theorem will enable us to prove discrete observability theorems for systems of vibrating strings and 
beams. For the latter, we will also need to investigate what happens when we add a new exponent to the system, i.e., when we consider sums of the form
\begin{equation}\label{16}
x(t)=x'e^{i\omega't}+\sum_{k=-\infty}^{\infty} x_ke^{i\omega_kt}
\end{equation}
with complex coefficients $x'$, $x_k$, instead of \eqref{11}, where $\omega'$ is a real number 
not belonging to the sequence $(\omega_k)$.

The following result is a discrete version of a generalization of a theorem of Haraux, allowing a weakened gap condition. In order to simplify its statement, let us introduce the quadratic form
\begin{equation*}
Q'(x):=\abs{x'}^2+Q(x).
\end{equation*}

\begin{theorem}\label{t3}
Assume \eqref{12} and introduce the sets $A_1$, $A_2$ as in Theorem \ref{t2}. 
Assume that for some positive $\delta>0$ and for some positive integer $J$ there exist two positive constants $c_1$, $c_2$, depending only on $\gam$ and $J\delta$, such that
\begin{equation}\label{17}
c_1Q(x)\le\delta\sum_{j=-J}^J\abs{x(j\delta)}^2\le c_2Q(x)
\end{equation}
for all sums of the form \eqref{11} with complex coefficients $x_k$ satisfying \eqref{15}. If $\omega'$ is a real number 
not belonging to the sequence $(\omega_k)$, then for every positive integer $J'$ there exist two positive constants $c_3$, 
$c_4$, depending only on $\gam$, $J\delta$, $J'\delta$ and
\begin{equation*}
\gam':=\inf_k\abs{\omega_k-\omega'}
\end{equation*}
and another constant $c'$, depending only on $\gam,$ and $J'\delta$, such that
\begin{equation}\label{18}
c_3Q'(x)\le\delta\sum_{j=-J-J'}^{J+J'}\abs{x(j\delta)}^2\le c_4Q'(x)
\end{equation}
for all sums of the form \eqref{16} with complex coefficients $x'$, $x_k$ satisfying \eqref{15} and
\begin{equation}\label{19}
\abs{\omega_k-\omega'}<2c'/\delta .
\end{equation} 
\end{theorem}

We may deduce from the preceding theorem the following

\begin{corollary}\label{c4}
Assume \eqref{12} and introduce the sets $A_1$, $A_2$ as in Theorem \ref{t2}. 
Assume that for $R>0$ there exist two positive constants $c_1$, $c_2$, depending only on $\gam$ and $R$, such that
\begin{equation*}
c_1Q(x)\le\int_{-R}^R \abs{x(t)}^2\ dt\le c_2Q(x)
\end{equation*}
for all sums of the form \eqref{11} with complex coefficients $x_k$ satisfying \eqref{15}. If $\omega'$ is a real number 
not belonging to the sequence $(\omega_k)$, then for every $R'>R$ there exist two positive constants $c_3$, 
$c_4$, depending only on $\gam$, $R$, $R'$ and
\begin{equation*}
\gam':=\inf_k\abs{\omega_k-\omega'}
\end{equation*}
and another constant $c'$, depending only on $\gam,$ and $R'$, such that
\begin{equation*}
c_3Q'(x)\le\int_{-R'}^{R'} \abs{x(t)}^2\ dt\le c_4Q'(x)
\end{equation*}
for all sums of the form \eqref{16} with complex coefficients $x'$, $x_k$ satisfying \eqref{15} and \eqref{19}.
\end{corollary}

Indeed, we may assume without loss of generality that $R'/R$ is a rational number. Then it suffices to apply Theorem \ref{t3} with arbitrarily large integers $J$ for which $J':=JR'/R$ is also integer, and with $\delta:=R/J$, and then letting $J\to\infty$.

Theorems \ref{t2} and \ref{t3} are proved in the following two sections. They are applied in Section \ref{s4} 
for the solution of two observability problems.

\section{Proof of Theorem \ref{t2}}\label{s2}

We proceed in three steps.

{\em First step.}  We begin by recalling the summatory formula of Poisson: if
$G$ is a
function belonging to $H_0^1(-\gam,\gam)$ and
its Fourier transform is given by the formula
\begin{equation*}
g(t)=\int_{-\infty}^{\infty}G(x) e^{-itx}\ dx
\end{equation*}
for all real $t$, then all functions of the form \eqref{11} with finitely many nonzero coefficients satisfy the following identity:
\begin{equation}\label{21}
\delta\sum_{j=-\infty}^{\infty}g(j\delta)\abs{x(j\delta)}^2=2\pi
\sum_{k,n=-\infty}^{\infty}G(\omega_k-\omega_n)x_k\overline{x_n}.
\end{equation}

For the proof we begin by remarking that since $\pi/\delta\ge\gam$, $G$ vanishes outside the interval
\begin{equation*}
I:=\Bigl(-\frac{\pi}{\delta},\frac{\pi}{\delta}\Bigr),
\end{equation*}
so that
\begin{equation*}
g(j\delta)=\int_IG(x) e^{-ij\delta x}\ dx
\end{equation*}
for all integers $j$. Since $G$ is H\"older continuous with exponent $1/2$, 
applying the Dini--Lipschitz theorem (see, e.g., \cite{Zyg}) to the trigonometric
orthonormal basis 
\begin{equation*}
\sqrt{\frac{\delta}{2\pi}}e^{ij\delta x},\quad j\in\ZZ
\end{equation*}
of $L^2(I)$, we conclude that
\begin{equation*}
\delta\sum_{j=-\infty}^{\infty}g(j\delta)e^{ij\delta x}=2\pi G_\delta (x)
\end{equation*}
for all real $x$, where $G_\delta$ denotes the $2\pi/\delta$-periodic
function which is equal to $G$ in the interval $I$. Observe that 
\begin{equation}\label{22}
 G_\delta(x)=G(x)\quad\text{whenever}\quad \abs{x}\le \frac{2\pi}{\delta}-\gam .
\end{equation}

Now we have
\begin{align*}
\delta\sum_{j=-\infty}^{\infty}g(j\delta)\abs{x(j\delta)}^2
&=\delta\sum_{k,n=-\infty}^{\infty}x_k\overline{x_n}\sum_{j=-\infty}^{\infty}g(j\delta)e^{i(\omega_k-\omega_n)j
\delta}\\
&=2\pi \sum_{k,n=-\infty}^{\infty}G_\delta(\omega_k-\omega_n)x_k\overline{x_n}\\
&=2\pi \sum_{k,n=-\infty}^{\infty}G(\omega_k-\omega_n)x_k\overline{x_n}.
\end{align*}
The last equality follows from \eqref{22} and from the fact that if $x_k\ne 0$ and $x_n\ne 0$, then by \eqref{15} we have necessarily 
\begin{equation*}
\abs{\omega_k-\omega_n}\le \frac{2\pi}{\delta}-\gam .
\end{equation*}
\medskip

{\em Second step.}  We prove the direct inequality (the second inequality in \eqref{14}. We are going to apply the identity 
\eqref{21} with the functions $H$, $G$ defined by
\begin{equation*}
H(x):=
\begin{cases}
\cos^2 \frac{\pi x}{2\gam}&\text{if $\abs{x}\le \gam$,}\\
0&\text{if $\abs{x}>\gam$,}
\end{cases}
\end{equation*}
the convolution product $G:=H*H$, and their Fourier transforms $h$ and $g$. One can readily verify (see \cite{BaiKomLor103} 
for details) that there exist two positive constants $\alpha$ and $\beta$ such that
\begin{align*}
&0\le G(0)-G(x)\le \alpha x^2\quad\text{for all}\quad x;\\
&G(x)=0\quad\text{whenever}\quad \abs{x}\ge\gam;\\
&g(t)\ge 0\quad\text{for all}\quad t;\\
&g(t)\ge \beta\quad\text{whenever}\quad \abs{t}\le\pi/(2\gam).
\end{align*}
We may assume without loss of generality that $\alpha\ge 1$.

Starting with \eqref{21} and using these relations we obtain the following estimates, where $J'$ denotes the (lower) integer part of 
$\pi/(2\gam\delta)$:
\begin{align*}
\frac{\beta}{2\pi}\delta\sum_{j=-J'}^{J'}\abs{x(j\delta)}^2
&\le \frac{\delta}{2\pi}\sum_{j=-\infty}^{\infty}g(j\delta)\abs{x(j\delta)}^2\\
&=\sum_{k,n=-\infty}^{\infty}G(\omega_k-\omega_n)x_k\overline{x_n}\\
&=\sum_{k\in A_1}G(0)\abs{x_k}^2+\sum_{k\in A_2}G(0)\bigl(\abs{x_k}^2+\abs{x_{k+1}}^2\bigr)\\
&\qquad +\sum_{k\in A_2}G(\omega_{k+1}-\omega_k)(x_k\overline{x_{k+1}}+\overline{x_k}x_{k+1})\\
&=\sum_{k\in A_1}G(0)\abs{x_k}^2+\sum_{k\in A_2}G(0)\abs{x_k+x_{k+1}}^2\\
&\qquad +\sum_{k\in A_2}(G(\omega_{k+1}-\omega_k)-G(0))(x_k\overline{x_{k+1}}+\overline{x_k}x_{k+1})\\
&\le \sum_{k\in A_1}G(0)\abs{x_k}^2+\sum_{k\in A_2}G(0)\abs{x_k+x_{k+1}}^2\\
&\qquad +\sum_{k\in A_2}(G(0)-G(\omega_{k+1}-\omega_k))\cdot (\abs{x_k}^2+\abs{x_{k+1}}^2\bigr)\\
&\le \sum_{k\in A_1}G(0)\abs{x_k}^2+\sum_{k\in A_2}G(0)\abs{x_k+x_{k+1}}^2\\
&\qquad +\sum_{k\in A_2}\alpha (\omega_{k+1}-\omega_k)^2\cdot (\abs{x_k}^2+\abs{x_{k+1}}^2\bigr)\\
&\le \alpha Q(x).
\end{align*}
We conclude that for $J=J'$ the direct inequality holds with
\begin{equation*}
c_2:=\frac{2\pi \alpha}{\beta} .
\end{equation*}

A translation argument in \cite{BaiKomLor103}, Remark 2.6 shows that we have, more generally,
\begin{equation*}
\delta\sum_{j=-J'}^{J'}\abs{x(t'+j\delta)}^2\le \frac{4\pi \alpha}{\beta}(1+\abs{t'}^2) Q(x)
\end{equation*}
for every real number $t'$.
The direct inequality for a general integer $J$ hence follows by covering the set $\set{-J,\ldots,J}$ of consecutive integers by $M$ 
translates of $\set{-J',\ldots,J'}$ where $M$ denotes the upper integer part of $(2J+1)/(2J'+1)$, and summing the  $M$ 
corresponding inequalities. 

{\em Third step.}  
For the proof of the {\em inverse inequality} let us introduce the same function $H$ as above, but define this time $G:=R^2H*H+H'*H'$. Denoting by $h$ and $g$ the 
Fourier transforms of $H$ and $G$, now we have the following properties with suitable positive constants $\alpha$ and 
$\beta$:
\begin{align*}
&G(0)-G(x)\ge \alpha x^2\quad\text{if}\quad \abs{x}\le\gam;\\
&G(x)=0\quad\text{whenever}\quad \abs{x}\ge\gam;\\
&G(0)>0\quad\text{and}\quad G(0)-G(x)>0\quad\text{for all}\quad x\ne 0;\\
&g(t)\le 0\quad\text{whenever}\quad \abs{t}\ge R;\\
&g(t)\le \beta\quad\text{for all}\quad t.
\end{align*}
We may assume without loss of generality that $\alpha\le G(0)$.

Applying \eqref{21} and using these relations we obtain the following estimates, where $J'$ denotes the upper integer part of 
$\pi/(2\gam\delta)$:
\begin{align*}
\frac{\beta}{2\pi}\delta\sum_{j=-J'}^{J'}\abs{x(j\delta)}^2
&\ge \frac{\delta}{2\pi}\sum_{j=-\infty}^{\infty}g(j\delta)\abs{x(j\delta)}^2\\
&=\sum_{k,n=-\infty}^{\infty}G(\omega_k-\omega_n)x_k\overline{x_n}\\
&=\sum_{k\in A_1}G(0)\abs{x_k}^2+\sum_{k\in A_2}G(0)\bigl(\abs{x_k}^2+\abs{x_{k+1}}^2\bigr)\\
&\qquad +\sum_{k\in A_2}G(\omega_{k+1}-\omega_k)(x_k\overline{x_{k+1}}+\overline{x_k}x_{k+1})\\
&=\sum_{k\in A_1}G(0)\abs{x_k}^2+\sum_{k\in A_2}G(0)\bigl(\abs{x_k}^2+\abs{x_{k+1}}^2\bigr)\\
&\qquad +\sum_{k\in A_2}G(\omega_{k+1}-\omega_k)\cdot 
(\abs{x_k+x_{k+1}}^2-\abs{x_k}^2-\abs{x_{k+1}}^2)\\
&=\sum_{k\in A_1}G(0)\abs{x_k}^2+\sum_{k\in A_2}(G(0)-G(\omega_{k+1}-\omega_k))\bigl(\abs{x_k}^2+\abs{x_{k+1}}^2\bigr)\\
&\qquad +\sum_{k\in A_2}G(\omega_{k+1}-\omega_k)\cdot \abs{x_k+x_{k+1}}^2.
\end{align*}
Putting $y:=\omega_{k+1}-\omega_k$ for brevity, it remains to show that
\begin{equation*}
\abs{x_k+x_{k+1}}^2+y^2\bigl(\abs{x_k}^2+\abs{x_{k+1}}^2\bigr)
\end{equation*}
is majorized by a constant multiple of
\begin{equation*}
G(y)\abs{x_k+x_{k+1}}^2+(G(0)-G(y))\bigl(\abs{x_k}^2+\abs{x_{k+1}}^2\bigr)
\end{equation*}
for all $0<y<\gam$. We show the stronger inequality
\begin{multline*}
\abs{x_k+x_{k+1}}^2+\frac{G(0)-G(y)}{\alpha }\bigl(\abs{x_k}^2+\abs{x_{k+1}}^2\bigr)\\
\le \frac{G(y)}{\alpha}\abs{x_k+x_{k+1}}^2+\frac{3(G(0)-G(y))}{\alpha}\bigl(\abs{x_k}^2+\abs{x_{k+1}}^2\bigr)
\end{multline*}
or equivalently, that
\begin{equation*}
\left( \alpha -G(y)\right)\abs{x_k+x_{k+1}}^2 \le 2 \left( G(0)-G(y)\right) \bigl(\abs{x_k}^2+\abs{x_{k+1}}^2\bigr).
\end{equation*}
This is obvious for $G(y)\ge \alpha$ because the right-hand side is nonnegative. If $G(y)< \alpha$, then the inequality follows from our assumption $\alpha\le G(0)$ and from the elementary estimate $\abs{x_k+x_{k+1}}^2\le 2\bigl(\abs{x_k}^2+\abs{x_{k+1}}^2\bigr)$.

\section{Proof of Theorem \ref{t3}}\label{s3}

\begin{proof}[\rm\bf Proof of the {\em direct} part of \eqref{18}] 
Applying the second inequality of \eqref{17} to the function
\begin{equation*}
z(t):= x(t)-x'e^{i\omega't}
\end{equation*}
instead of $x(t)$, we obtain that
\begin{align*}
\delta\sum_{j=-J}^{J}\abs{x(j\delta)}^2 
&\le 2\delta\sum_{j=-J}^{J}\abs{z(j\delta)}^2 +2\delta\sum_{j=-J}^{J}\bigl\lvert x'e^{i\omega'j\delta}\bigr\rvert^2\\
&\le 2c_2Q(z)+2\delta (2J+1)\abs{x'}^2 \\
&\le \max\set{2c_2,6J\delta} Q'(x).
\end{align*}
Using \cite{BaiKomLor103}, Remark 2.6 this inequality implies that, more generally,
\begin{equation}\label{31}
\delta\sum_{j=m-J}^{m+J}\abs{x(j\delta)}^2 \le \max\set{4c_2,12J\delta}\bigl(1+\abs{m\delta}^2\bigr) Q'(x)
\end{equation}
for every integer $m$. (In order to use this remark, we also apply Remark 2.5 of that paper which enables us to choose 
$0<\gam_0\le\gam$ sufficiently small so that $\abs{\omega'-\omega_k}<\gam_0$ for all $k$. Then $\omega'$ belongs to $A_1$ in the extended exponent set, so that the coresponding quadratic form is $Q'(x)$.)

Now the second inequality of \eqref{18} follows easily by covering the set $\set{-J-J',\ldots, J+J'}$ of consecutive integers by $M:=\lceil (2J+2J'+1)/(2J+1) \rceil$ translates of $\set{-J,\ldots, J}$, and summing the corresponding inequalities \eqref{31}. Since $\abs{m}\le J'$ in all these inequalities, we obtain that
\begin{equation*}
\delta\sum_{j=-J-J'}^{J+J'}\abs{x(j\delta)}^2 \le  c_4Q'(x)
\end{equation*}
with
\begin{equation*}
c_4=\Bigl(1+\frac{2J+2J'+1}{2J+1}\Bigr)
\max\set{4c_2,12J\delta}\bigl(1+\abs{J'\delta}^2\bigr).\qedhere
\end{equation*}
\end{proof}

\begin{proof}[\rm\bf Proof of the {\em inverse} part of \eqref{18}] 
For $x$ given by \eqref{16}, the formula
\begin{equation*}
y(t):=x(t)-\frac{1}{2J'}\sum_{n=-J'}^{J'-1}e^{-i\omega' n\delta}x(t+n\delta)
\end{equation*}
defines a function $y$ of the form \eqref{11}: an easy computation shows that
\begin{equation*}
y(t)=\sum_{k=-\infty}^{\infty }\Bigl[1-\frac{1}{2J'}\sum_{n=-J'}^{J'-1}e^{i(\omega_k-\omega')n\delta} \Bigr]x_ke^{i\omega_kt}
=:\sum_{k=-\infty}^{\infty }y_ke^{i\omega_kt}.
\end{equation*}

Observe that
\begin{align*}
\frac{1}{2J'}\sum_{n=-J'}^{J'-1}e^{i(\omega_k-\omega')n\delta} 
&=\frac{1}{2J'}\frac{e^{i(\omega_k-\omega')J'\delta}-e^{-i(\omega_k-\omega')J'\delta}}
{e^{i(\omega_k-\omega')\delta}-1}\\
&=\frac{\sin (\omega_k-\omega')J'\delta}{J'\bigl(e^{i(\omega_k-\omega')\delta}-1\bigr)}
\end{align*}
and therefore
\begin{align*}
\Bigl\vert\frac{1}{2J'}\sum_{n=-J'}^{J'-1}e^{i(\omega_k-\omega')n\delta} \Bigr\vert
&=\Bigl\vert\frac{\sin (\omega-\omega')J'\delta}{J'\bigl(e^{i(\omega-\omega')\delta}-1\bigr)}
\Bigr\vert\\
&=\Bigl\vert\frac{\sin(\omega_k-\omega')J'\delta}
{2J'\sin (\omega_k-\omega')\delta/2}\Bigr\vert\\
&=\Bigl\vert\frac{\sin(\omega_k-\omega')J'\delta}{(\omega_k-\omega')J'\delta}\Bigr\vert
\cdot \Bigl\vert\frac{(\omega_k-\omega')\delta/2}{\sin(\omega_k-\omega')\delta/2}\Bigr\vert
=:\eps_k.
\end{align*}

For the sequel we need the following

\begin{lemma}\label{l5}\mbox{}

(a) There exists a constant $c'$, depending only on $\inf_k\abs{\omega_k-\omega'}$ and $J'\delta$, such that
\begin{equation*}
\eps:=\sup_k\eps_k<1
\end{equation*}
where the supremum is taken over the indices $k$ satisfying \eqref{19}.
\smallskip

(b) The function
\begin{equation*}
\omega\mapsto f(\omega):=\frac{\sin (\omega-\omega')J'\delta}{J'\bigl(e^{i(\omega-\omega')\delta}-1\bigr)}
\end{equation*}
is Lipschitzian in the interval $(\omega'-2c'/\delta,\omega'+2c'/\delta)$ with some constant $L$, depending only on $J'\delta$.
\end{lemma}

\begin{proof}\mbox{}

(a) Since
\begin{equation*}
\inf_k\abs{\omega_k-\omega'}>0,
\end{equation*}
we have
\begin{equation*}
\eps'=\eps'(J'\delta):=\sup_k\Bigl\vert
\frac{\sin(\omega_k-\omega')J'\delta}{(\omega_k-\omega')J'\delta}
\Bigr\vert<1.
\end{equation*}
It suffices to choose $c'>0$ sufficiently small so that
\begin{equation*}
\inf_{0<\abs{x}<c'}\Bigl\vert\frac{\sin x}{x}\Bigr\vert>\eps'.
\end{equation*}

(b) First we note that under the condition \eqref{19} we have
\begin{equation*}
\Bigl\vert\frac{e^{i(\omega-\omega')\delta}-1}{\delta}\Bigr\vert
=\Bigl\vert\frac{\sin(\omega-\omega')\delta/2}{\delta/2}\Bigr\vert
>\eps'\abs{\omega-\omega'}\ge \eps'\gam'.
\end{equation*}
Therefore
\begin{align*}
\abs{f'(\omega )}
&=\Bigl\vert\frac{J'\delta \cos (\omega-\omega')J'\delta}{J'\bigl(e^{i(\omega-\omega')\delta}-1\bigr)}
-\frac{\bigl[\sin (\omega-\omega')J'\delta\bigr]iJ'\delta e^{i(\omega-\omega')\delta}}
{\bigl[J'\bigl(e^{i(\omega-\omega')\delta}-1\bigr)\bigr]^2}\Bigr\vert\\
&\le \Bigl\vert\frac{\delta}{e^{i(\omega-\omega')\delta}-1}\Bigr\vert
+\frac{1}{J'\delta}\Bigl\vert\frac{\delta}{e^{i(\omega-\omega')\delta}-1}\Bigr\vert^2\\
&\le \frac{1}{\eps'\gam'}+\frac{1}{(J'\delta)(\eps'\gam')^2}.
\end{align*}
The Lipschitz property follows by the mean value theorem because the constant on the right-hand depends only on $\gam'$ and $J'\delta$.
\end{proof}

It follows from part (a) of the lemma that
\begin{equation}\label{32}
\abs{x_k-y_k}\le \eps\abs{x_k}\quad\text{for all $k$ satisfying \eqref{19}}.
\end{equation}

We claim that
\begin{equation}\label{33}
Q(x)\le c_0Q(y)
\end{equation}
with a suitable constant $c_0$ (depending on $\eps$). If $k\in A_1$, then we deduce from \eqref{32} that
\begin{equation}\label{34}
\abs{x_k}^2\le (1-\eps)^{-2}\abs{y_k}^2.
\end{equation}
using part (b) of the above, in case $k\in A_2$ we have
\begin{align*}
\abs{(x_k-y_k)+(x_{k+1}-y_{k+1})}
&\le\abs{f(\omega_k)}\cdot \abs{x_k+x_{k+1}} 
+\abs{f(\omega_k)-f(\omega_{k+1})}\cdot \abs{x_{k+1}}\\
&\le \eps \abs{x_k+x_{k+1}}+L\abs{\omega_k-\omega_{k+1}}\cdot\abs{x_{k+1}}.
\end{align*}
Hence
\begin{equation*}
(1-\eps)\abs{x_k+x_{k+1}}
\le \abs{y_k+y_{k+1}}+L\abs{\omega_k-\omega_{k+1}}\cdot(1-\eps)^{-1}\abs{y_{k+1}}
\end{equation*}
and therefore
\begin{equation*}
\abs{x_k+x_{k+1}}
\le (1-\eps)^{-1}\abs{y_k+y_{k+1}}+L\abs{\omega_k-\omega_{k+1}}\cdot(1-\eps)^{-2}\abs{y_{k+1}}.
\end{equation*}
Using this relation we obtain that
\begin{multline}\label{35}
\abs{x_k+x_{k+1}}^2+(\omega_k-\omega_{k+1})^2\bigl(\abs{x_k}^2+\abs{x_{k+1}}^2\bigr)\\
\le
2(1-\eps)^{-2}\abs{y_k+y_{k+1}}^2\\
+\bigl(2L^2(1-\eps)^{-4}+(1-\eps)^{-2}\bigr)
\cdot (\omega_k-\omega_{k+1})^2\bigl(\abs{y_k}^2+\abs{y_{k+1}}^2\bigr).
\end{multline}
Finally, \eqref{33} follows from \eqref{34} and \eqref{35}.

Next we show that
\begin{equation}\label{36}
\sum_{j=-J}^{J}\abs{y(j\delta)}^2\le 4\sum_{m=-J-J'}^{J+J'-1}\abs{x(m\delta)}^2.
\end{equation} 
Indeed, using the Cauchy--Schwarz inequality we have
\begin{align*}
\abs{y(t)}^2
&\le 2\abs{x(t)}^2+2\Bigl\vert\frac{1}{2J'}\sum_{n=-J'}^{J'-1}e^{-i\omega' n\delta}x(t+n\delta)\Bigr\vert^2\\
&\le 2\abs{x(t)}^2+\frac{2}{4(J')^2}\Bigl(\sum_{n=-J'}^{J'-1}\abs{e^{-i\omega' n\delta}}^2\Bigr)
\cdot\Bigl(\sum_{n=-J'}^{J'-1}\abs{x(t+n\delta)}^2\Bigr)\\
&= 2\abs{x(t)}^2+\frac{1}{J'}\sum_{n=-J'}^{J'-1}\abs{x(t+n\delta)}^2
\end{align*}
for every $t$. Hence
\begin{align*}
\sum_{j=-J}^{J}\abs{y(j\delta)}^2
&\le 2\sum_{j=-J}^{J}\abs{x(j\delta)}^2
+ \frac{1}{J'}\sum_{j=-J}^{J}\sum_{n=-J'}^{J'-1}\abs{x(j\delta+n\delta)}^2\\
&\le 2\sum_{j=-J}^{J}\abs{x(j\delta)}^2
+ 2\sum_{m=-J-J'}^{J+J'-1}\abs{x(m\delta)}^2\\
&\le 4\sum_{m=-J-J'}^{J+J'-1}\abs{x(m\delta)}^2.
\end{align*}

Now applying the first inequality of \eqref{17} for $y$ instead of $x$ and using \eqref{33} and \eqref{36} we obtain that
\begin{equation}\label{37}
Q(x)\le c_0Q(y)\le \frac{c_0}{c_1}\delta\sum_{j=-J}^{J}\abs{y(j\delta)}^2
\le \frac{4c_0}{c_1}\delta\sum_{m=-J-J'}^{J+J'-1}\abs{x(m\delta)}^2. 
\end{equation} 
Furthermore, using the function $z$ as introduced in the proof of the direct inequality,
\begin{align}
\abs{x'}^2
&=\frac{1}{2J+2J'}\sum_{m=-J-J'}^{J+J'-1}\bigl\vert x'e^{i\omega' m\delta}\bigr\vert^2\notag\\
&\le \frac{1}{J+J'}\sum_{m=-J-J'}^{J+J'-1}\abs{x(m\delta )}^2+\abs{z(m\delta )}^2.\label{38}
\end{align}
Applying to $z$ the already proved direct inequality and then the inequality  \eqref{37}, we obtain that
\begin{equation*}
\delta\sum_{m=-J-J'}^{J+J'-1}\abs{z(m\delta )}^2
\le c_4Q'(z)
= c_4Q(x)
\le \frac{4c_0c_4}{c_1}\delta\sum_{m=-J-J'}^{J+J'-1}\abs{x(m\delta )}^2.
\end{equation*}
Combining this with \eqref{38} we get
\begin{equation}\label{39}
\abs{x'}^2\le \Bigl(1+\frac{4c_0c_4}{c_1}\Bigr)\frac{1}{J+J'}\sum_{m=-J-J'}^{J+J'-1}\abs{x(m\delta )}^2.
\end{equation} 
Finally, we conclude from \eqref{37} and \eqref{39} that
\begin{equation*}
Q'(x)\le c\delta \sum_{m=-J-J'}^{J+J'-1}\abs{x(m\delta )}^2
\end{equation*}
avec
\begin{equation*}
c:=\max\Bigl\lbrace
\frac{4c_0}{c_1},\Bigl(1+\frac{4c_0c_4}{c_1}\Bigr)\frac{1}{J\delta +J'\delta }
\Bigr\rbrace .\qedhere
\end{equation*}
\end{proof}

\section{Simultaneous observability of strings and beams}\label{s4}

Fix a number $0<a<1$ arbitrarily and consider the following problem:
\begin{equation}
\begin{cases}
u_{tt}-u_{xx}=0\quad \text{in}\quad (0,a)\times {\RR},\\
u_{tt}-u_{xx}=0\quad \text{in}\quad (a,1)\times {\RR},\\
u(0,\cdot)=u(a,\cdot)=u(1,\cdot)=0\quad \text{in}\quad {\RR},\\
u(\cdot,0)=u_{0,a}\quad \text{and}\quad u_t(\cdot,0)=u_{1,a}\quad
\text{in}\quad (0,a),\\
u(\cdot,0)=u_{0,1-a}\quad \text{and}\quad u_t(\cdot,0)=u_{1,1-a}\quad
\text{in}\quad (a,1).
\end{cases}
\label{41}
\end{equation}
This problem is well posed for 
\begin{equation*}
u_{0,a}\in H_0^1(0,a),\quad u_{1,a}\in L^2(0,a)\quad u_{0,1-a}\in H_0^1(a,1),\quad 
\text{and}\quad u_{1,1-a}\in L^2(a,1).
\end{equation*}
The following result follows at once from the proof of part (a) Theorem 5.1 in \cite{BaiKomLor103} if we apply Theorem 
\ref{t2} above instead of its continuous version. Set
\begin{equation*}
\gamma=\frac{\pi}{2}\min\Bigl\{\frac{1}{a},\frac{1}{1-a}\Bigr\}.
\end{equation*}

\begin{theorem}\label{t6}
The following estimate holds for almost every $0<a<1$ and for every $\eps>0$. Given $0<\delta\le \pi/\gam$ arbitrarily, fix 
an integer $J$ such that  $J\delta>2\max\set{a,1-a}$. For every $t'\in\RR$, the inequality
\begin{equation*}
\norm{u_0}_{H^{-\eps}(0,1)}^2+
\norm{u_1}_{H^{-1-\eps}(0,1)}^2
\le C\delta\sum_{j=-J}^{J}\abs{u_x(a-0,t'+j\delta)-u_x(a+0,t'+j\delta)}^2
\end{equation*}
is satisfied for all solutions of \eqref{41} whose initial data are $u_{0,a}$, $u_{1,a}$ linear combinations of the basis 
functions
\begin{equation*}
\sin(n\pi a^{-1}x),\quad n\le \frac{a}{\delta}-\frac{1}{4}\min\Bigl\{1,\frac{a}{1-a}\Bigr\}
\end{equation*}
and whose initial data are $u_{0,1-a}$, $u_{1,1-a}$ linear combinations of the basis functions
\begin{equation*}
\sin(m\pi(1-a)^{-1}x),\quad  m\le \frac{1-a}{\delta}-\frac{1}{4}\min\Bigl\{1,\frac{1-a}{a}\Bigr\}.
\end{equation*}
\end{theorem}

Now consider the following problem:
\begin{equation}
\begin{cases}
u_{tt}+u_{xxxx}=0\quad \text{in}\quad (0,a)\times{\RR},\\
u_{tt}+u_{xxxx}=0\quad \text{in}\quad (a,1)\times{\RR},\\
u(0,\cdot)=u(a,\cdot)=u(1,\cdot)=0\quad \text{in}\quad {\RR},\\
u_{xx}(0,\cdot)=u_{xx}(a,\cdot)=u_{xx}(1,\cdot)=0\quad \text{in}\quad {\RR},\\
u(\cdot,0)=u_{0,a}\quad \text{and}\quad u_t(\cdot,0)=u_{1,a}\quad
\text{in}\quad (0,a),\\
u(\cdot,0)=u_{0,1-a}\quad \text{and}\quad u_t(\cdot,0)=u_{1,1-a}\quad
\text{in}\quad (a,1).
\end{cases}
\label{42}
\end{equation}
This system models two vibrating beams with simply supported endpoints, one of which is common to both 
beams. It is well posed for initial data satisfying
\begin{equation*}
u_{0,a}\in H_0^1(0,a),\quad u_{1,a}\in H^{-1}(0,a)\quad u_{0,1-a}\in H_0^1(a,1)\quad 
\text{and}\quad u_{1,1-a}\in H^{-1}(a,1).
\end{equation*}
The following result follows at once from the proof of part (a) Theorem 6.1 in \cite{BaiKomLor103} if we apply Theorem 
\ref{t2} above instead of its continuous version:

\begin{theorem}\label{t7}
The following estimate holds for almost every $0<a<1$ and for every $\eps>0$. Given $\gam>0$ and $0<\delta\le \pi/\gam$ 
arbitrarily, fix an integer $J$ such that  $J\delta>\pi/\gam$. For every $t'\in\RR$, the inequality
\begin{equation*}
\norm{u_0}_{H^{1-\eps}(0,1)}^2+
\norm{u_1}_{H^{-1-\eps}(0,1)}^2
\le C\delta\sum_{j=-J}^{J}\abs{u_x(a-0,t'+j\delta)-u_x(a+0,t'+j\delta)}^2
\end{equation*}
is satisfied for all solutions of \eqref{42} whose initial data are $u_{0,a}$, $u_{1,a}$ linear combinations of the basis 
functions
\begin{equation*}
\sin(n\pi a^{-1}x),\quad n\le \frac{a}{\pi}\sqrt{\frac{\pi}{\delta}-\frac{\gam}{2}}
\end{equation*}
and whose initial data are $u_{0,1-a}$, $u_{1,1-a}$ linear combinations of the basis functions
\begin{equation*}
\sin(m\pi(1-a)^{-1}x),\quad  m\le \frac{1-a}{\pi}\sqrt{\frac{\pi}{\delta}-\frac{\gam}{2}}.
\end{equation*}
\end{theorem}

\end{document}